%% file: HSZ11a.tex
\documentclass{svmult}

\usepackage{amsmath}
\usepackage{amssymb}
\usepackage{amsfonts}
\usepackage{amsbsy}
\usepackage{amscd}
\usepackage{amstext}
\usepackage{mathptmx}
\usepackage{helvet}
\usepackage{courier}
\usepackage{bbm}
\usepackage{bm}
\usepackage{dsfont}
\usepackage{stmaryrd}
\usepackage{latexsym}
\usepackage[english]{babel}
\usepackage{mathrsfs}
\usepackage[numbers]{natbib}
\usepackage{chapterbib}
\usepackage{color}
\usepackage{makeidx}
\usepackage{graphicx}
\usepackage{graphics}
\usepackage{epsfig}
\usepackage{psfrag}
\usepackage{subfigure}
\usepackage{wrapfig}
\usepackage{multicol}
\usepackage{multirow}
\usepackage[bottom]{footmisc}
\usepackage{url}
\usepackage{array}
\usepackage{slashbox}
\usepackage{pict2e}
\usepackage{enumerate}
\usepackage[pagewise]{lineno} 
\usepackage{verbatim}
\usepackage{oldgerm}
\usepackage{hangcaption}
\usepackage[ruled]{algorithm}
\usepackage{algorithmic}
\usepackage{threeparttable}
\usepackage{longtable}
\usepackage{epic}
\usepackage{epsf}
\usepackage[latin1]{inputenc}
\usepackage[T1]{fontenc}
\usepackage{setspace}
\usepackage{varioref}
\usepackage{textcomp}
\usepackage{fancyhdr}
\usepackage{dcolumn}
\usepackage{booktabs}
\usepackage{listings}
\usepackage{enumitem}

  \newcounter{mnote}
  \let\oldmarginpar\marginpar
    \renewcommand\marginpar[1]{\-\oldmarginpar[\raggedleft\footnotesize #1]%
    {\raggedright\footnotesize #1}}

\newcommand{\tbarXszy}{|\hspace*{-0.15em}|\hspace*{-0.15em}|}

\begin{document}
\input{HSZ11a_body}

\end{document}

%% file: HSZ11a_body.tex

\title*{Adaptive Finite Element Methods with Inexact Solvers for the Nonlinear Poisson-Boltzmann Equation}
\titlerunning{Inexact AFEM for Nonlinear PBE}

\author{Michael Holst\inst{1}\and
Ryan Szypowski\inst{2}\and
Yunrong Zhu\inst{3}}

\institute{
 Departments of Mathematics and Physics, University of California San Diego, 
     La Jolla, CA 92093. 
 Supported in part by NSF Awards~0715146 and 0915220,
     and by CTBP and NBCR,
     \texttt{mholst@math.ucsd.edu},
     \texttt{http://ccom.ucsd.edu/\~{}mholst/}
\and
 Department of Mathematics, University of California San Diego, 
     La Jolla, CA 92093. 
 Supported in part by NSF Award~0715146, \texttt{rszypows@math.ucsd.edu}
\and
 Department of Mathematics, University of California San Diego, 
     La Jolla, CA 92093.
 Supported in part by NSF Award~0715146, \texttt{zhu@math.ucsd.edu}
}

\maketitle


\section{Introduction}

In this article we study adaptive finite element methods (AFEM) with inexact
solvers for a class of semilinear elliptic interface problems.
We are particularly interested in nonlinear problems with discontinuous 
diffusion coefficients, such as the nonlinear Poisson-Boltzmann equation 
and its regularizations. The algorithm we study consists of the standard
SOLVE-ESTIMATE-MARK-REFINE procedure common to many adaptive finite element
algorithms, but where the SOLVE step involves only a full solve on the 
coarsest level, and the remaining levels involve only single Newton updates 
to the previous approximate solution.
We summarize a recently developed AFEM convergence theory for inexact 
solvers appearing in~\cite{BHSZ11b}, and present a sequence of numerical 
experiments that give evidence that the theory does in fact predict the
contraction properties of AFEM with inexact solvers.
The various routines used are all designed to maintain a linear-time
computational complexity.

An outline of the paper is as follows.
In Section~\ref{sec:1-szypowski}, we give a brief overview of the
Poisson-Boltzmann equation.
In Section~\ref{sec:2-szypowski}, we describe AFEM algorithms, and
introduce a variation involving inexact solvers.
In Section~\ref{sec:3-szypowski}, we give a sequence of numerical experiments
that support the theoretical statements on convergence and optimality.
Finally, in Section~\ref{sec:4-szypowski} we make some final observations.

\section{Regularized Poisson-Boltzmann Equation}
\label{sec:1-szypowski}

We use standard notation for Sobolev spaces. In particular, we denote $\|\cdot\|_{0,G}$ the $L^{2}$ norm on any subset $G\subset \mathbb{R}^{3},$ and denote $\|\cdot \|_{1,2, G}$ the $H^{1}$ norm on $G$. 

\begin{figure}[t]
\begin{center}
\includegraphics[width=1.8in]{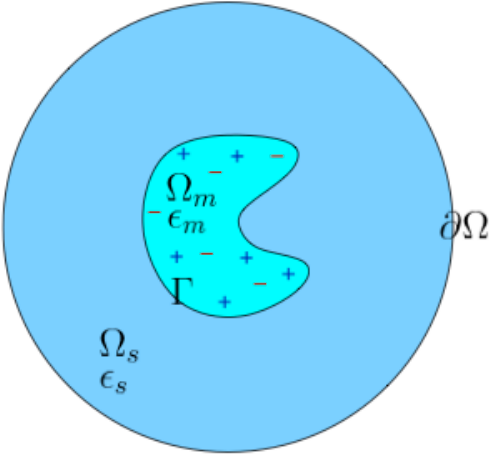}
\caption{Schematic of a molecular domain.}
\end{center}
\label{fig:mol}
\end{figure}

Let $\Omega:=\Omega_{m}\cup \Gamma \cup \Omega_{s}$ be a bounded Lipschitz domain in $\mathbb{R}^{3}$, which consists of the molecular region $\Omega_{m},$ the solvent region $\Omega_{s}$ and their interface $\Gamma:= \overline{\Omega}_{m}\cap\overline{\Omega}_{s}$ (see Figure~\ref{fig:mol}).
Our interest in this paper is to solve the following regularized Poisson-Boltzmann equation 
in the weak form: find $u \in  H_g^1(\Omega):=\{u\in H^{1}(\Omega): u|_{\partial \Omega}=g\}$ such that
\begin{equation}
	\label{eqn:pbe_weak}
	a(u, v) + ( b(u), v ) = (f, v)
    \quad \forall v\in H_0^1(\Omega),
\end{equation}
where
$a(u,v) = \int_{\Omega}\epsilon \nabla u \cdot\nabla v dx,$ $(b(u),v ) =   \int_{\Omega}\kappa^2
\sinh(u)v dx$. Here we assume that the diffusion coefficient
$\epsilon$ is piecewise positive constant $\epsilon|_{\Omega_{m}} =
\epsilon_{m}$ and $\epsilon|_{\Omega_{s}} = \epsilon_{s}$.  The
modified Debye-H\"uckel parameter $\kappa^{2}$ is also piecewise
constant with $\kappa^{2}(x)|_{\Omega_{m}} = 0$ and
$\kappa^{2}(x)|_{\Omega_{s}} >0$.  The equation \eqref{eqn:pbe_weak} arises from several regularization schemes (cf. \cite{Chen.L;Holst.M;Xu.J2007,CLW03}) of the nonlinear Poisson-Boltzmann equation:
\begin{equation*}
-\nabla \cdot (\epsilon \nabla u) + \kappa^2 \sinh u = \sum_{i=1}^{N} z_i \delta(x_i),
\end{equation*}
where the right hand side represents $N$ fixed points with charges
$z_{i}$ at positions $x_{i}$, and $\delta$ is the Dirac delta
distribution.  

It is easy to verify that the bilinear form in~\eqref{eqn:pbe_weak} 
satisfies:
\begin{equation*}
c_{0} \|u\|_{1,2}^2 \leq a(u,u),
\qquad
a(u,v) \leq c_{1} \|u\|_{1,2} \|v\|_{1,2}, 
\qquad \forall u,v \in H_0^{1}(\Omega),
\label{eqn:coercive-bounded}
\end{equation*}
where $0 < c_{0} \leq c_{1} < \infty$ are constants depending only on $\epsilon$.
These properties imply the norm on $H_0^{1}(\Omega)$ is equivalent
to the energy norm 
$\tbarXszy\cdot\tbarXszy \colon H_0^{1}(\Omega) \rightarrow \mathbb{R}$,
\begin{equation*}
	\tbarXszy u \tbarXszy^2 = a(u,u),
	\qquad
	c_{0} \|u\|_{1,2}^2 \leq \tbarXszy u \tbarXszy^2 \leq c_{1} \|u\|_{1,2}^2.
   \label{eqn:equiv}
\end{equation*}

Let $\mathcal{T}_{h}$ be a shape-regular conforming triangulation of $\Omega$, and let $V_{g}(\mathcal{T}_{h}):=\{v\in H^{1}(\Omega): v|_{\tau} \in \mathbb{P}_{1}(\tau)\;\; \forall \tau \in \mathcal{T}_{h}\}$ be the standard piecewise linear finite element space defined on $\mathcal{T}_{h}$. For simplicity, we assume that the interface $\Gamma$ is resolved by $\mathcal{T}_{h}$. Then the finite element approximation of \eqref{eqn:pbe_weak} reads: 
find $u_h \in
V_{g}(\mathcal{T}_{h})$ such that
\begin{equation}
\label{eqn:fem}
a(u_h, v) + (b(u_h), v) = (f, v), \quad \forall v\in V_{0}(\mathcal{T}_{h}).
\end{equation}

We close this section with a summary of a priori $L^{\infty}$ bounds for the solution $u$ to \eqref{eqn:pbe_weak} and the discrete solution $u_{h}$ to \eqref{eqn:fem}, which play a key role in the finite element error analysis of \eqref{eqn:fem} and adaptive algorithms. For interested reader, we refer to \cite{Chen.L;Holst.M;Xu.J2007,Holst.M;McCammon.J;Yu.Z;Zhou.Y2009} for details.
\begin{theorem} 
\label{thm:Linfty}
There exist $u_{+}, u_{-}\in L^{\infty}(\Omega)$ such that the solution $u$ of
\eqref{eqn:pbe_weak} satisfies the following a priori $L^{\infty}$ bounds:
\begin{equation}
\label{eqn:cont-Linfty}
	u_{-} \le u \le u_{+},  \qquad \mbox { a.e.  in } \Omega. 
\end{equation}
Moreover, if the triangulation $\mathcal{T}_{h}$ satisfies that 
\begin{equation}
	\label{eqn:mesh}
	a(\phi_{i}, \phi_{j}) \le -\frac{\sigma}{h^{2}} \sum_{e_{i,j} \subset \tau} |\tau|, \quad \mbox{for some}\quad \sigma >0,
\end{equation}
for all the adjacent vertices $i\neq j$ with the basis function $\phi_{i}$ and $\phi_{j},$ then the discrete solution $u_{h}$ of \eqref{eqn:fem} also has the a priori $L^{\infty}$ bound
\begin{equation}
\label{eqn:dis-Linfty}
	 \| u_{h}\|_{L^{\infty}(\Omega)} \le C,
\end{equation}
where $C$ is a constant independent of $h.$
\end{theorem}
We note that the mesh condition is generally not needed practically,
and in fact can also be avoided in analysis for certain
nonlinearites~\cite{BHSZ11b}.

\section{Adaptive FEM with Inexact Solvers}
\label{sec:2-szypowski}

Given a discrete solution $u_{h}\in V_{g}(\mathcal{T}_{h})$, let us define the residual based error indicator $\eta(u_{h}, \tau)$: 
$$
	\eta^{2}(u_{h}, \tau)= h_{\tau}^{2} \|b(u_{h}) -f\|_{0,\tau}^{2} + \sum_{e\subset \partial \tau} h_{e}\|[(\epsilon \nabla u_{h})\cdot n_{e}]\|^{2}_{0,e},
$$
where $[(\epsilon \nabla u_{h})\cdot n_{e}]$ denote the jump of the flux across a face $e$ of $
\tau.$ For any subset $\mathcal{S}\subset \mathcal{T}_{h},$ we set $\eta^{2}(u_{h}, \mathcal{S}) :=\sum_{\tau\in \mathcal{S}} \eta^{2}(u_{h}, \tau).$ By using the a priori $L^{\infty}$ bounds Theorem ~\ref{thm:Linfty}, we can show (cf. \cite{Holst.M;McCammon.J;Yu.Z;Zhou.Y2009}) that the error indicator satisfies:
\begin{equation}
\label{eqn:upper}
\tbarXszy u - u_h\tbarXszy^2 \le C_{1} \eta^2(u_h, \hat{\mathcal{T}}_h);
\end{equation}
and 
\begin{equation}
\label{eqn:lipschitz}
|\eta(v, \tau) - \eta(w, \tau)| \le C_{2}\tbarXszy v - w \tbarXszy_{\omega_\tau}, \quad \forall v,w\in V_{g}(\mathcal{T}_{h})
\end{equation}
 where $\omega_\tau = \cup_{\tau'\in \mathcal{T}_{h}, \bar{\tau}'\cap \bar{\tau} \neq \emptyset} \tau'$ and $\tbarXszy v\tbarXszy_{\omega_\tau}^{2} = \int_{\omega_{\tau}} \epsilon |\nabla v|^{2} dx$.

Given an initial triangulation $\mathcal{T}_{0}$, the standard adaptive finite element method (AFEM) generates a sequence $\left[u_{k}, \mathcal{T}_{k}, \{\eta(u_{k}, \tau)\}_{\tau\in \mathcal{T}_{k}}\right] $ based on the iteration of the form:
$$
\hbox{$\textsf{SOLVE}$} \rightarrow \hbox{$\textsf{ESTIMATE}$} \rightarrow \hbox{$\textsf{MARK}$} \rightarrow \hbox{$\textsf{REFINE}$}.
$$
Here the \textsf{SOLVE} subroutine is usually assumed to be exact, namely $u_{k}$ is the exact solution to the nonlinear equation \eqref{eqn:fem}; the \textsf{ESTIMATE} routine computes the element-wise residual indicator $\eta(u_{k}, \tau)$; the \textsf{MARK} routine uses standard D\"{o}rfler marking (cf. \cite{Dorfler.W1996}) where
$\mathcal{M}_k \subset \mathcal{T}_k$ is chosen so that
$$
\eta(u_k, \mathcal{M}_k) \ge \theta \eta(u_k, \mathcal{T}_k)
$$
for some parameter $\theta\in (0, 1];$ finally, the routine \textsf{REFINE} subdivide the marked elements and possibly some neighboring elements in certain way such that the new triangulation preserves shape-regularity and conformity. 

During last decade, a lot of theoretical work has been done to show the convergence of the AFEM with exact solver (see \cite{Nochetto.R;Siebert.K;Veeser.A2009} and the references cited therein for linear PDE case, and \cite{Holst.M;Tsogtgerel.G;Zhu.Y2010} for nonlinear PDE case). To the best of the authors knowledge, there are only a couple of convergence results of AFEM for symmetric linear elliptic equations (cf. \cite{Stevenson.R2007,Arioli.M;Georgoulis.E;Loghin.D2009}) which take the numerical error into account. To distinct with the exact solver case, we use $\hat{u}_{k}$ and $\hat{\mathcal{T}}_{k}$ to denote the numerical approximation to \eqref{eqn:fem} and the triangulation obtained from the adaptive refinement using the inexact solutions. 

Due to the page limitation, we only state the main convergence result of the AFEM with inexact solver for solving \eqref{eqn:pbe_weak} below. More detailed analysis and extension are reported in \cite{BHSZ11b}.
\begin{theorem}
\label{thm:main}
Let $\{\hat{\mathcal{T}}_k, \hat{u}_k\}_{k \ge 0}$ be the sequence of
meshes and approximate solutions computed by the AFEM algorithm. Let $u$ denote
the exact solution and $u_k$ denote the exact discrete solutions on
the meshes $\hat{\mathcal{T}}_k$. Then, there exist constants $\mu > 0$, $\nu \in (0, 1)$, $\gamma > 0$,
and $\alpha \in (0,1)$ such that if the inexact solutions satisfy
\begin{equation}
\label{eqn:approx_prop}
\mu \tbarXszy u_k - \hat{u}_k \tbarXszy^2 +  \tbarXszy u_{k+1} - \hat{u}_{k+1} \tbarXszy^2 \le \nu
\eta^2(\hat{u}_k, \hat{\mathcal{T}}_k)
\end{equation}
then
\begin{equation}
\label{eqn:contraction}
 \tbarXszy u - u_{k + 1} \tbarXszy^2 + \gamma \eta^2(\hat{u}_{k+1},
\hat{\mathcal{T}}_{k+1})\le \alpha^2(\tbarXszy u - u_k\tbarXszy^2 + \gamma \eta^2(\hat{u}_k,
\hat{\mathcal{T}}_k)).
\end{equation}
Consequently, $\lim_{k\to \infty} u_{k} = \lim_{k\to \infty} \hat{u}_{k} = u.$
\end{theorem}
The proof of this theorem is based on the upper bound \eqref{eqn:upper} of the exact solution, the Lipschitz property \eqref{eqn:lipschitz} of the error indicator, D\"orfler marking, and the following quasi-orthogonality between the exact solutions:
\begin{equation}
\label{eqn:quasi-orthogonality}
\tbarXszy u - u_{k+1}\tbarXszy^2 \le \Lambda \tbarXszy u - u_k \tbarXszy^2 - \tbarXszy u_{k + 1} - u_k \tbarXszy^2
\end{equation}
where $\Lambda$ can be made close to 1 by refinement. For a proof of the inequality \eqref{eqn:quasi-orthogonality}, see for example \cite{Holst.M;McCammon.J;Yu.Z;Zhou.Y2009}.

To achieve the optimal computational complexity, we should avoid solving the nonlinear system \eqref{eqn:fem} as much as we could. The two-grid algorithm \cite{Xu.J1996b} shows that a nonlinear solver on a coarse grid combined with a Newton update on the fine grid still yield quasi-optimal approximation. Motivated by this idea, we propose the following AFEM algorithm with inexact solver, which contains only one nonlinear solver on the coarsest grid, and Newton updates on each follow-up steps: 
\begin{algorithm}
\caption{:~$\left[\hat{u}_{k}, \hat{\mathcal{T}}_{k}, \{\eta(\hat{u}_k, \tau)\}_{\tau\in\hat{\mathcal{T}_k}}\right] :=$ \textsf{Inexact\_AFEM}$(\mathcal{T}_{0}, \theta)$}
\label{alg:iafem}

\raggedright


\textbf{\scriptsize 1~~} $\hat{u}_{0}=u_{0} :=$ \textsf{NSOLVE}$(\mathcal{T}_{0})$  $\quad$ \%{\it Nonlinear solver on initial triangulation}\;

\textbf{\scriptsize 2~~} \textbf{for} $k := 0, 1, \cdots$ \textbf{do}
{

\textbf{\scriptsize 3~~} \hspace*{0.4cm} $\{\eta(\hat{u}_k, \tau)\}_{\tau\in\hat{\mathcal{T}_k}} :=
\textsf{ESTIMATE}(\hat{u}_k, \hat{\mathcal{T}}_k)$\;

\textbf{\scriptsize 4~~} \hspace*{0.4cm} $\mathcal{M}_k := \textsf{MARK}(\{\eta(\hat{u}_k,
\tau)\}_{\tau\in\hat{\mathcal{T}_k}}, \hat{\mathcal{T}}_k, \theta)$\;

\textbf{\scriptsize 5~~} \hspace*{0.4cm} $\hat{\mathcal{T}}_{k+1} := \textsf{REFINE}(\hat{\mathcal{T}}_k,
\mathcal{M}_k)$\;

\textbf{\scriptsize 6~~} \hspace*{0.4cm} $\hat{u}_{k+1} := \textsf{UPDATE}(\hat{u}_{k}, \hat{\mathcal{T}}_{k+1})$ $\quad$\%{\it One-step Newton update}\;  
}

\textbf{\scriptsize 7~~} \textbf{end}

\end{algorithm}
In Algorithm~\ref{alg:iafem}, the \textsf{NSOLVE} routine is used only on the coarsest mesh and is
implemented using Newton's method run to certain convergence tolerance.  For the rest of the solutions, a single step of Newton's method is used to update the previous approximation.  That is,
\textsf{UPDATE} computes $\hat{u}_{k+1}$ such that
\begin{equation}
\label{eqn:update}
a(\hat{u}_{k+1} - \hat{u}_{k}, \phi) + (b'(\hat{u}_{k})(\hat{u}_{k+1} -
\hat{u}_{k}), \phi) = 0
\end{equation}
for every $\phi \in V(\hat{\mathcal{T}}_{k+1})$.
We remark that since \eqref{eqn:update} is only a linear problem, we could use the local multilevel method to solve it in (near) optimal complexity (cf. \cite{Chen.L;Holst.M;Xu.J;Zhu.Y2010}). Therefore, the overall computational complexity of the Algorithm~\ref{alg:iafem} is nearly optimal. 

We should point out that it is not obvioius how to enforce the required approximation property~\eqref{eqn:approx_prop} that $\hat{u}_{k}$ must satisfy for the theorem.
This is examined in more detail in~\cite{BHSZ11b}.
However, numerical evidence in the following section shows Algorithm~\ref{alg:iafem} is an efficient algorithm, and the results matches the ones from AFEM with exact solver. 
%

%
%

\section{Numerical Experiments}
\label{sec:3-szypowski}

In this section we present some numerical experiments to illustrate the 
result in Theorem~\ref{thm:main}, implemented with FETK.
The software utilizes the standard piecewise-linear finite element
space for discretizing \eqref{eqn:pbe_weak}.  Algorithm~\ref{alg:iafem}
is implemented with care taken to guarantee that each of the steps
runs in linear time relative to the number of vertices in the mesh.
The linear solver used is Conjugate Gradients preconditioned by
diagonal scaling.  The estimator is computed using a high-order
quadrature rule, and, as mentioned above, the marking strategy is
D\"{o}rfler marking where the estimated error have been binned to
maintain linear complexity while still marking the elements with the
largest error.  Finally, the refinement is longest edge bisection,
with refinement outside of the marked set to maintain conformity of
the mesh.

We present three sets of results in order to explore the effects of
the inexact solver in multiple contexts.  For all problems, we present
a convergence plot using both inexact and exact solvers (including a
reference line of order $N^{-\frac{1}{3}}$) as well as a
representative cut-away of a mesh with around 30,000 vertices.  The
exact discrete solution is computed using the standard AFEM algorithm
where the solution on each mesh is computed by allowing Newton's
method to continue running to convergence to within a tolerance of
$10^{-8}$.  This modifies not only the solution on a given mesh, but
also the sequence of meshes generated, since the algorithm may mark
different simplices.  However, in all cases, convergence is identical
to high precision and meshes place the unknowns as expected.

\begin{figure}[tbh]
\begin{center}
\includegraphics[width=1.8in]{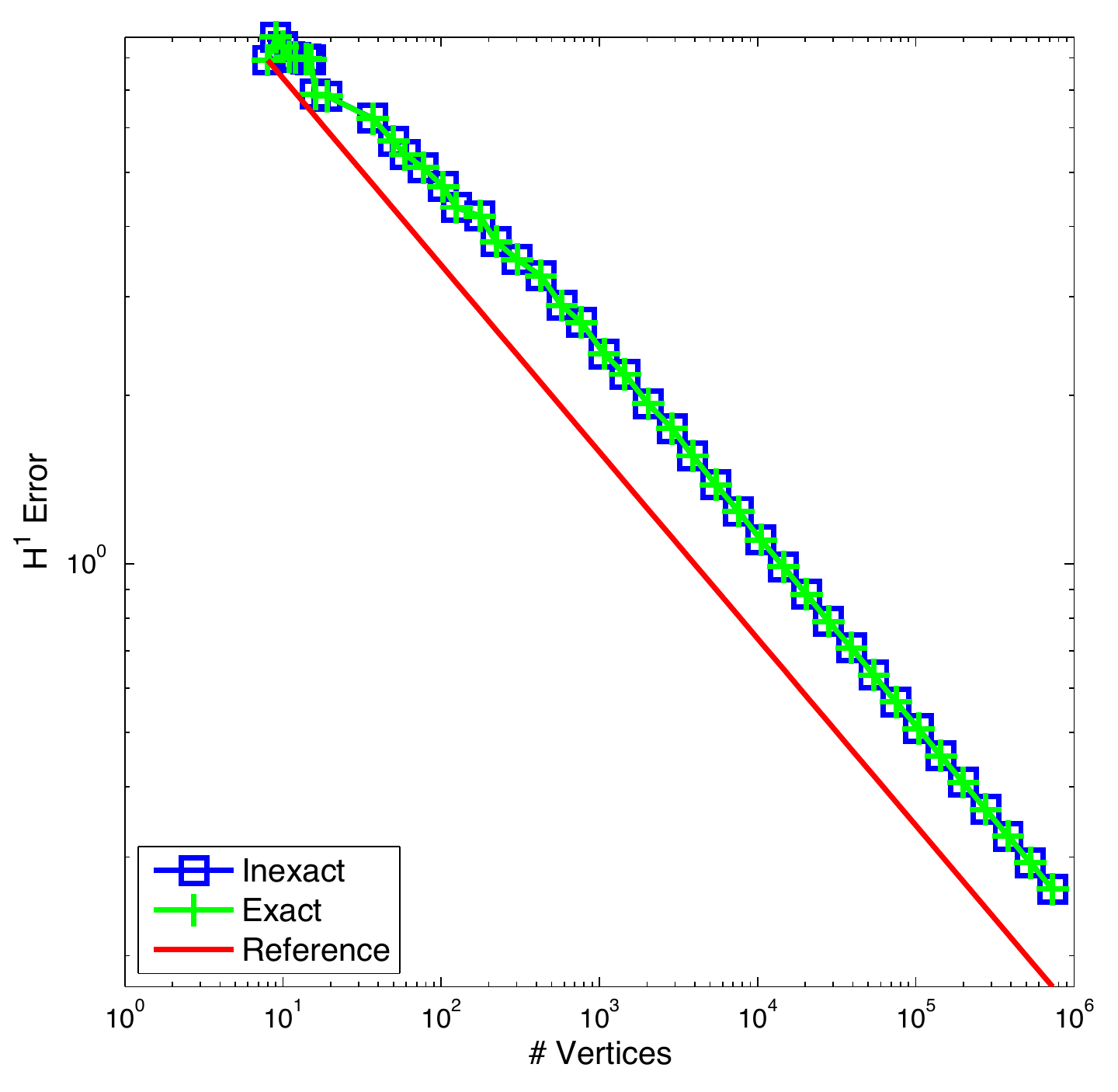}\includegraphics[width=2.7in]{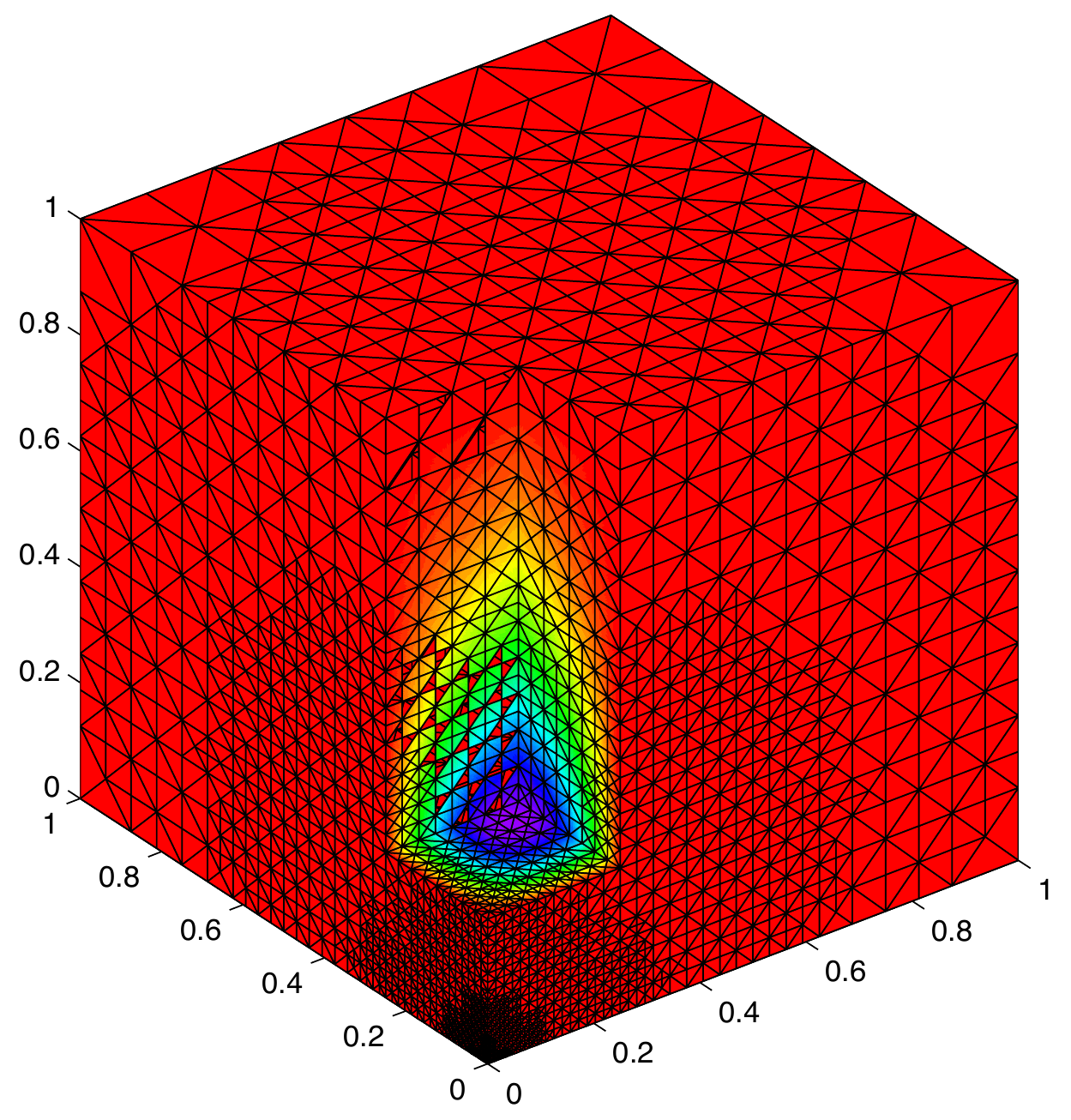}
\end{center}
\caption{\label{fig:corner}Convergence plot and mesh cut-away for the
  corner singularity problem.}
\end{figure}

The first result uses constant coefficients across the entire domain
$\Omega = [0,1]^3$ and we choose a right hand side so that the
derivative of the exact solution is large near the origin. The
boundary conditions chosen for this problem are homoegenous Dirichlet
boundary conditions. Specifically, the exact solution is given by $u =
u_1 u_2$ where
$$
u_1 = \sin(\pi x)\sin(\pi y)\sin(\pi z)
$$
is chosen to satisfy the boundary condition and
$$
u_2 = (x^2 + y^2 + x^2 + 10^{-4})^{-1.5}.
$$
The results can be seen in Figure~\ref{fig:corner}.

\begin{figure}
\begin{center}
\includegraphics[width=1.8in]{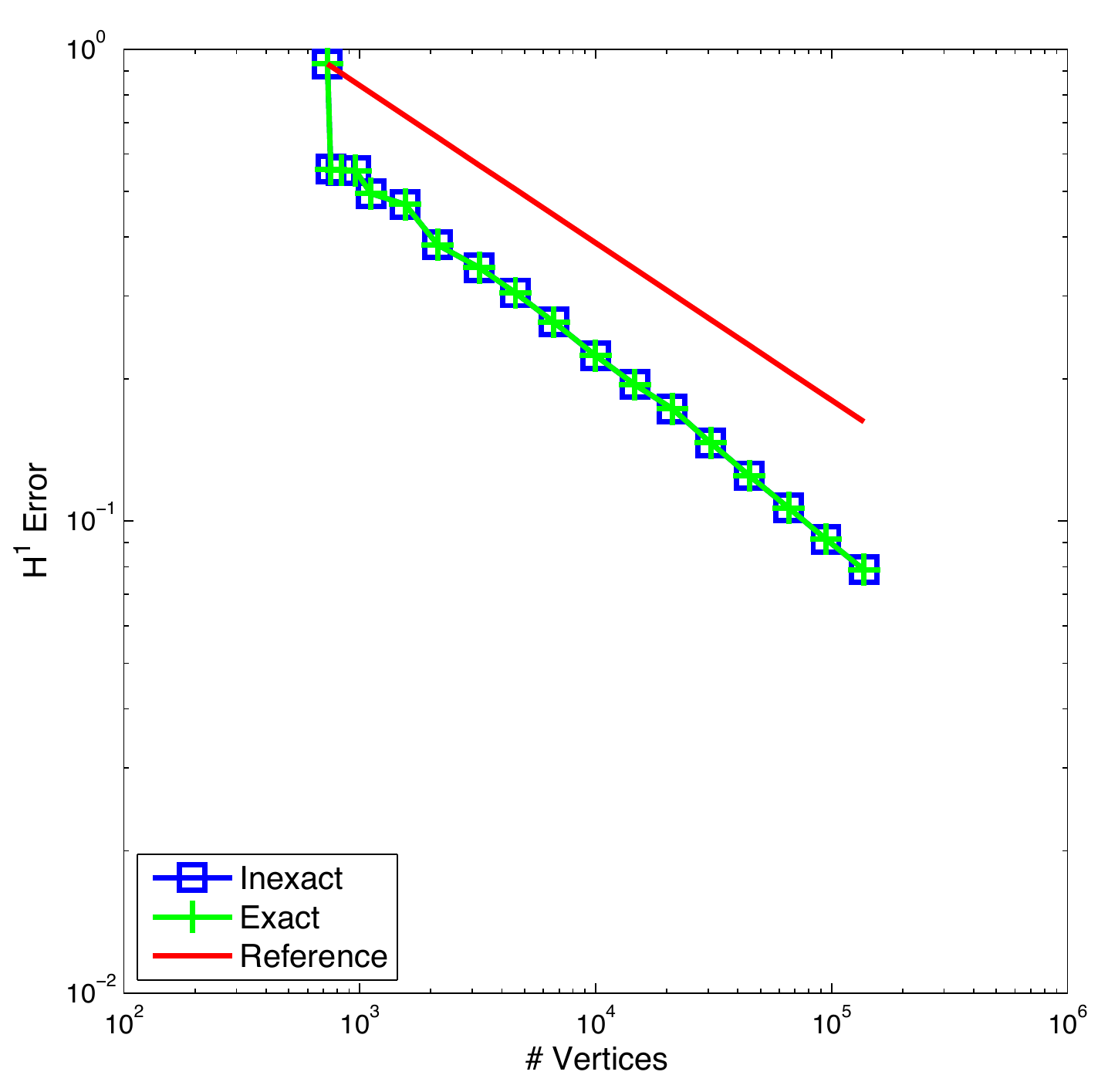}\includegraphics[width=2.7in]{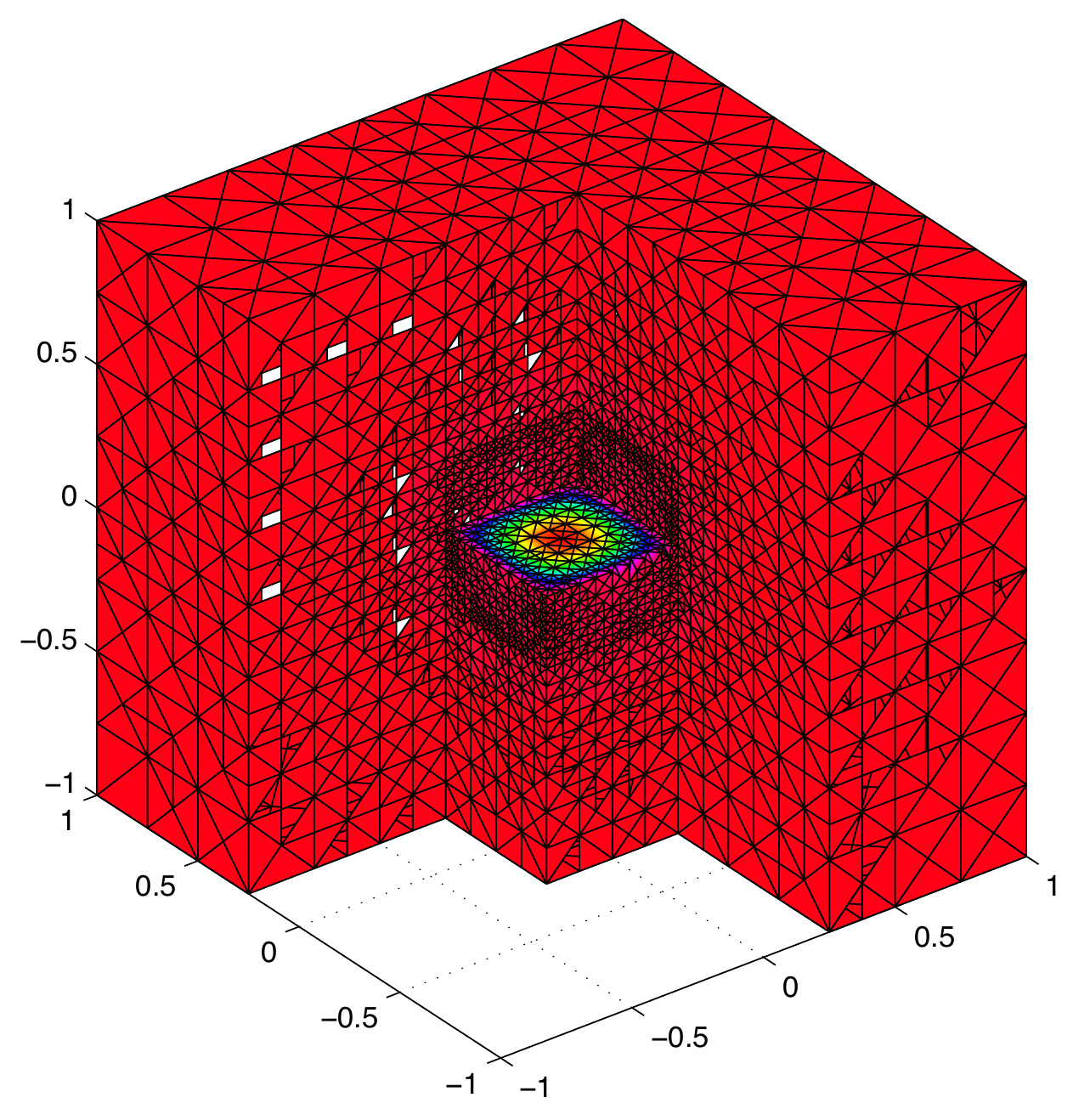}
\end{center}
\caption{\label{fig:rpbe1}Convergence plot and mesh cut-away for the
  Poisson-Boltzmann problem.}
\end{figure}

The second result uses the domain $\Omega = [-1,1]^3$ and $\Omega_m =
\left[-\frac{1}{4}, \frac{1}{4}\right]$ with constants $\epsilon_s = 80
, \epsilon_m = 2, \kappa_s =1,$ and $\kappa_m = 0$. Homogeneous Neumann
conditions are chosen for the boundary and the right hand side is
simplified to a constant.  Because an exact solution is unavailable
for this (and the following) problem, the error is computed by
comparing to a discrete solution on a mesh with around 10 times the
number of vertices as the finest mesh used in the adaptive algorithm.
Figure~\ref{fig:rpbe1} shows the results for this problem.  As can be
seen the refinement favors the interface and the inexact and exact
solvers perform as expected.

\begin{figure}
\begin{center}
\includegraphics[width=1.8in]{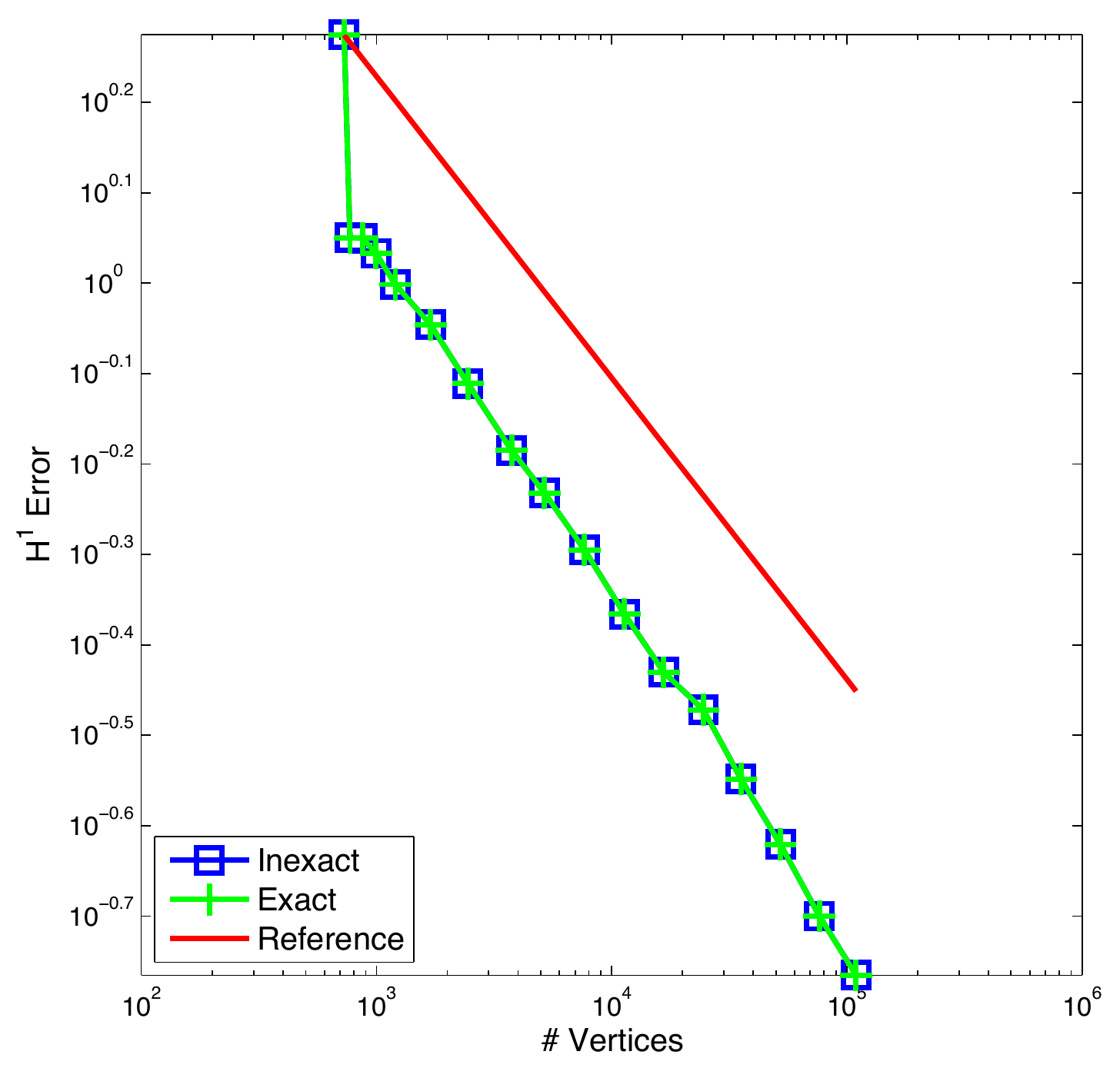}\includegraphics[width=2.7in]{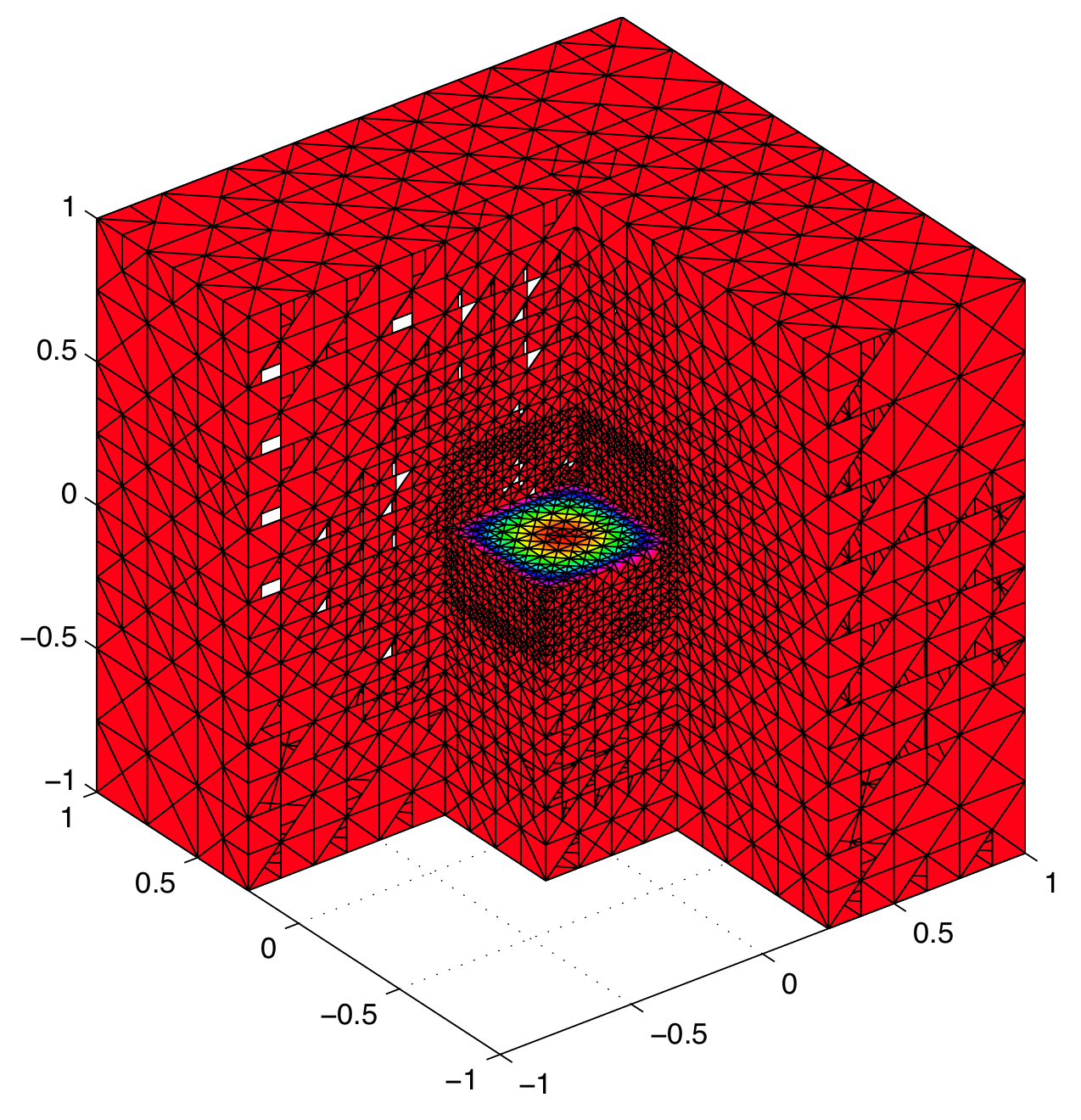}
\end{center}
\caption{\label{fig:rpbe2}Convergence plot and mesh cut-away for the
 exaggerated-jump Poisson-Boltzmann problem.}
\end{figure}

The final result is chosen to test the robustness of the inexact
method to large coefficient jumps.  The domain and boundary conditions
are the same as the previous example, but the constants chosen as are
$\epsilon_s = 1000, \epsilon_m = 10, \kappa_s =1,$ and $\kappa_m =
0$.  The results can be seen in Figure~\ref{fig:rpbe2}, and they
show a scenario very similar to the previous example, with the
refinement restricted even further to the interface.

\section{Conclusion}
\label{sec:4-szypowski}

In this article we have studied AFEM with inexact solvers for a class of 
semilinear elliptic interface problems with discontinuous diffusion 
coefficients, with emphasis on the nonlinear Poisson-Boltzmann equation.
The algorithm we studied consisted of the standard
SOLVE-ESTIMATE-MARK-REFINE procedure common to many adaptive finite element
algorithms, but where the SOLVE step involves only a full solve on the 
coarsest level, and the remaining levels involve only single Newton updates 
to the previous approximate solution.
The various routines used are all designed to maintain a linear-time 
computational complexity. 
Our numerical results indicate that the recently developed AFEM convergence 
theory for inexact solvers in~\cite{BHSZ11b} does predict the actual
behavior of the methods.


\bibliographystyle{plainnat}
\bibliography{HSZ11a}
